# Метод гармонического баланса для отыскания приближённых периодических решений системы Лоренца

**Александр Николаевич ПЧЕЛИНЦЕВ**[1],
**Андрей Андреевич ПОЛУНОВСКИЙ**[2], **Ирина Юрьевна ЮХАНОВА**[1]

[1] ФГБОУ ВО «Тамбовский государственный технический университет»
392000, Российская Федерация, г. Тамбов, ул. Советская, 106
ORCID: https://orcid.org/0000-0003-4136-1227, e-mail: pchelintsev.an@yandex.ru
ORCID: https://orcid.org/0000-0002-8339-0459, e-mail: irina_yu_10@mail.ru

[2] ФГБОУ ВО «Московский государственный технический университет им. Н.Э. Баумана
(национальный исследовательский университет)»
105005, Российская Федерация, г. Москва, 2-я Бауманская ул., 5
ORCID: https://orcid.org/0000-0002-6557-3649, e-mail: apap2009@yandex.ru

# The harmonic balance method for finding approximate periodic solutions of the Lorenz system

**Alexander N. PCHELINTSEV**[1],
**Andrey A. POLUNOVSKIY**[2], **Irina Yu. YUKHANOVA**[1]

[1] Tambov State Technical University
106 Sovetskaya St., Tambov 392000, Russian Federation
ORCID: https://orcid.org/0000-0003-4136-1227, e-mail: pchelintsev.an@yandex.ru
ORCID: https://orcid.org/0000-0002-8339-0459, e-mail: irina_yu_10@mail.ru

[2] Bauman Moscow State Technical University
5 Baumanskaya 2-ya, Moscow 105005, Russian Federation
ORCID: https://orcid.org/0000-0002-6557-3649, e-mail: apap2009@yandex.ru

**Аннотация.** Рассматривается метод гармонического баланса для нахождения приближённых периодических решений динамической системы Лоренца. При разработке программного обеспечения, реализующего описываемый метод, был выбран математический пакет Maxima. Показаны недостатки символьных вычислений для получения системы нелинейных алгебраических уравнений относительно циклической частоты, постоянных членов и амплитуд гармоник, составляющих искомое решение. Для ускорения расчётов впервые эта система была получена в общем виде. Приведены результаты вычислительного эксперимента – коэффициенты тригонометрических полиномов, приближающих найденное периодическое решение, начальное условие и период цикла. Полученные результаты были проверены с помощью описанного ранее в работах авторов высокоточного метода интегрирования, основанного на аппроксимации степенными рядами.

**Abstract.** We consider the harmonic balance method for finding approximate periodic solutions of the Lorenz system. When developing software that implements the described method, the math package Maxima was chosen. The drawbacks of symbolic calculations for obtaining a system of nonlinear algebraic equations with respect to the cyclic frequency, free terms and amplitudes of the harmonics, that make up the desired solution, are shown. To speed up the calculations, this system was obtained in a general form for the first time. The results of the computational experiment are given: the coefficients of trigonometric polynomials approximating the found periodic solution, the initial condition, and the cycle period. The results obtained were verified using a high-precision method of numerical integration based on the power series method and described earlier in the articles of the authors.

## Введение

Рассмотрим нелинейную систему дифференциальных уравнений, введённую Э. Лоренцом в работе [1],

$$\begin{cases} \dot{x}_1 = \sigma(x_2 - x_1), \\ \dot{x}_2 = rx_1 - x_2 - x_1 x_3, \\ \dot{x}_3 = x_1 x_2 - bx_3, \end{cases} \quad (0.1)$$

где $\sigma = 10$, $r = 28$, $b = 8/3$ (эти значения параметров системы (0.1) теперь называют классическими, а саму систему — системой Лоренца).

В статье [1] для системы Лоренца с классическими параметрами доказано следующее утверждение: существует такое $C > 0$, что для любого решения

$$X(t) = [x_1(t)\ x_2(t)\ x_3(t)]^{\mathrm{T}}$$

начиная с некоторого момента времени становится справедливым неравенство $|X(t)| \leq C$, и дивергенция векторного поля скоростей системы (0.1) отрицательна всюду в $\mathbb{R}^3$. Тогда [1] существует предельное множество — аттрактор Лоренца, — к которому притягиваются все траектории динамической системы при $t \to \infty$. Таким образом, аттрактор определяет поведение решений динамической системы на больших отрезках времени.

У. Такер в работе [2] доказал гиперболичность аттрактора в системе (0.1), т. е. аттрактор состоит из траекторий, всюду плотных на нём (континуум седловых циклов), вдоль которых близкие траектории экспоненциально разбегаются; это и создает их хаотическое поведение. Тогда (как отмечает Д.В. Аносов в послесловии к книге [3, с. 285]) *в аттракторе системы* (0.1) *может существовать бесконечное число асимптотически устойчивых периодических траекторий*, но их область притяжения может быть достаточно малой (трудно улавливаемой в численном эксперименте).

Как известно (см., например, [4, 5]), символическую динамику используют для отслеживания циклов в системе Лоренца. Разбивают область в фазовом пространстве, содержащую аттрактор, на конечное число подобластей. Обозначая каждый элемент разбиения буквой, траектории на аттракторе, проходящие через соответствующие области, кодируются последовательностями таких символов. Если в последовательности имеется регулярность — повторяемость групп символов, — то соответствующая траектория считается циклом. Однако возвращаемость траектории в некоторую окрестность своей части не говорит о её замкнутости. Критику результатов подобных вычислительных экспериментов можно найти, например, в [6].

В 2004 г. Д. Вишванат опубликовал работу [7], в которой привел начальные условия и периоды для трёх циклов в аттракторе Лоренца с достаточно большой точностью. Алгоритм вычислений основан на методе Линдштедта–Пуанкаре (ЛП), на который (в отличие от методов численного интегрирования) не влияет устойчивость цикла, к которому строятся приближения. Полученные в [7] вычислительные данные можно проверить, решая задачу Коши высокоточными численными методами (см, например, [9]).

Анализ работ [7, 8] Д. Вишваната показал, что автор приводит общее описание алгоритма без ссылок на программную реализацию (в MATLAB, как указано в его статьях). При этом не ясно, как для ЛП-метода символьно решается получаемая неоднородная линейная система дифференциальных уравнений с периодическими коэффициентами (например, для уравнения Ван дер Поля это сделать можно без особых проблем). Таким образом, актуальной задачей остаётся разработка алгоритма поиска циклов системы (0.1), детальное описание его реализации, получение начальных значений и периода цикла с заданной точностью.

Целью данной работы является отыскание приближённых периодических решений в системе Лоренца на основе метода гармонического баланса, являющегося более простым в реализации, чем ЛП-метод, используемый в [7, 8]. При этом будет получена в общем виде система нелинейных алгебраических уравнений относительно циклической частоты, постоянных членов и амплитуд гармоник, составляющих искомое решение.

## 1. Метод гармонического баланса

Попытки построить приближённые периодические решения системы (0.1) предпринимались и до Д. Вишваната (см., например, [10]) методом гармонического баланса, но с малой точностью представления вещественных чисел, при этом в статье [10] не указаны начальные условия и периоды найденных циклов (приведены только рисунки с цикла-ми). Сейчас этот метод активно развивается в работах А. Луо [11–13] для отыскания периодических решений нелинейных систем дифференциальных уравнений.

Будем использовать метод гармонического баланса для получения приближений к периодическим решениям системы (0.1). Для этого сделаем аппроксимацию фазовых координат на периоде $T$ тригонометрическими полиномами в общем виде с неизвестной циклической частотой $\omega$ (поскольку мы не знаем значение $T$; в общем случае оно может быть иррациональным числом):

$$x_1(t) \approx \tilde{x}_1(t) = x_{1,0} + \sum_{i=1}^{h} \left( c_{1,i} \cos(i\omega t) + s_{1,i} \sin(i\omega t) \right),$$

$$x_2(t) \approx \tilde{x}_2(t) = x_{2,0} + \sum_{i=1}^{h} \left( c_{2,i} \cos(i\omega t) + s_{2,i} \sin(i\omega t) \right),$$

$$x_3(t) \approx \tilde{x}_3(t) = x_{3,0} + \sum_{i=1}^{h} \left( c_{3,i} \cos(i\omega t) + s_{3,i} \sin(i\omega t) \right),$$

где $h$ — заданное количество гармоник. Если $i > h$, то мы полагаем

$$c_{1,i} = s_{1,i} = c_{2,i} = s_{2,i} = c_{3,i} = s_{3,i} = 0. \qquad (1.2)$$

В силу правой части системы (0.1) составим невязки

$$\begin{aligned}
\delta_1(t) &= \tilde{x}'_1(t) - \sigma[\tilde{x}_2(t) - \tilde{x}_1(t)], \\
\delta_2(t) &= \tilde{x}'_2(t) - [r\tilde{x}_1(t) - \tilde{x}_2(t) - \tilde{x}_1(t)\tilde{x}_3(t)], \\
\delta_3(t) &= \tilde{x}'_3(t) - [\tilde{x}_1(t)\tilde{x}_2(t) - b\tilde{x}_3(t)],
\end{aligned}$$

где штрихом переобозначена производная функции по времени. Если производить вычисления в аналитическом виде, то для каждой невязки нужно следующее:

1. Продифференцировать по времени соответствующий тригонометрический полином.
2. Где имеются произведения фазовых координат, перемножить соответствующие тригонометрические полиномы, преобразовав при этом произведения тригонометрических функций в суммы.
3. Привести подобные слагаемые для каждой функции cos() и sin() с соответствующим аргументом.
4. В силу равенств (1.2), отсечь от полученной невязки гармоники более высокого порядка.
5. Приравнять полученную невязку к нулю, т. е. коэффициенты при её гармониках.

Если собрать в единое целое найденные алгебраические уравнения для каждой невязки, то получим пока ещё незамкнутую систему нелинейных уравнений относительно неизвестных амплитуд $c_{1,i}$, $s_{1,i}$, $c_{2,i}$, $s_{2,i}$, $c_{3,i}$ и $s_{3,i}$ ($i = \overline{1,h}$), постоянных членов $x_{1,0}$, $x_{2,0}$ и $x_{3,0}$ и циклической частоты $\omega$. Количество неизвестных в системе равно $3(1+2h)+1 = 6h+4$, а уравнений — на единицу меньше.

Дополнительное уравнение можно получить исходя из следующих соображений. Известно (см. [5,7]), что искомые циклы пересекают плоскость, проходящую через положения равновесия системы (0.1)

$$O_1\left(-\sqrt{b(r-1)},\, -\sqrt{b(r-1)},\, r-1\right),\ O_2\left(\sqrt{b(r-1)},\, \sqrt{b(r-1)},\, r-1\right) \qquad (1.3)$$

и параллельную плоскости $x_1 O x_2$ (сечение Пуанкаре). Таким образом, третья координата в начальном условии для искомых циклов равна величине $r-1$, откуда $\tilde{x}_3(0) = r-1$. Тогда дополнительное уравнение системы имеет вид:

$$x_{3,0} + \sum_{i=1}^{h} c_{3,i} - 27 = 0.$$

Других дополнительных сведений о периодических решениях системы Лоренца авторы не встречали. Заметим, что для трёх циклов, найденных Д. Вишванатом, в начальном условии для третьей координаты было взято число 27.

Далее приведем пример системы уравнений при $h=2$:

$$\begin{cases}
\omega s_{1,1} - 10 c_{2,1} + 10 c_{1,1} = 0, \\
-10 s_{2,1} + 10 s_{1,1} - c_{1,1}\omega = 0, \\
2\omega s_{1,2} - 10 c_{2,2} + 10 c_{1,2} = 0, \\
-10 s_{2,2} + 10 s_{1,2} - 2 c_{1,2}\omega = 0, \\
10 x_{1,0} - 10 x_{2,0} = 0, \\
c_{1,1} x_{3,0} + c_{3,1} x_{1,0} + \dfrac{s_{1,1} s_{3,2}}{2} + \dfrac{s_{1,2} s_{3,1}}{2} + \omega s_{2,1} + \dfrac{c_{1,1} c_{3,2}}{2} + \dfrac{c_{1,2} c_{3,1}}{2} + c_{2,1} - 28 c_{1,1} = 0, \\
s_{1,1} x_{3,0} + s_{3,1} x_{1,0} + \dfrac{c_{1,1} s_{3,2}}{2} - \dfrac{c_{1,2} s_{3,1}}{2} + s_{2,1} + \dfrac{c_{3,1} s_{1,2}}{2} - \dfrac{c_{3,2} s_{1,1}}{2} - 28 s_{1,1} - c_{2,1}\omega = 0, \\
c_{1,2} x_{3,0} + c_{3,2} x_{1,0} - \dfrac{s_{1,1} s_{3,1}}{2} + 2\omega s_{2,2} + \dfrac{c_{1,1} c_{3,1}}{2} + c_{2,2} - 28 c_{1,2} = 0, \\
s_{1,2} x_{3,0} + s_{3,2} x_{1,0} + \dfrac{c_{1,1} s_{3,1}}{2} + s_{2,2} - 28 s_{1,2} + \dfrac{c_{3,1} s_{1,1}}{2} - 2 c_{2,2}\omega = 0, \\
x_{1,0} x_{3,0} + x_{2,0} - 28 x_{1,0} + \dfrac{s_{1,2} s_{3,2}}{2} + \dfrac{s_{1,1} s_{3,1}}{2} + \dfrac{c_{1,2} c_{3,2}}{2} + \dfrac{c_{1,1} c_{3,1}}{2} = 0, \\
-c_{1,1} x_{2,0} - c_{2,1} x_{1,0} + \omega s_{3,1} - \dfrac{s_{1,1} s_{2,2}}{2} - \dfrac{s_{1,2} s_{2,1}}{2} + \dfrac{8 c_{3,1}}{3} - \dfrac{c_{1,1} c_{2,2}}{2} - \dfrac{c_{1,2} c_{2,1}}{2} = 0, \\
-s_{1,1} x_{2,0} - s_{2,1} x_{1,0} + \dfrac{8 s_{3,1}}{3} - \dfrac{c_{1,1} s_{2,2}}{2} + \dfrac{c_{1,2} s_{2,1}}{2} - \dfrac{c_{2,1} s_{1,2}}{2} + \dfrac{c_{2,2} s_{1,1}}{2} - c_{3,1}\omega = 0, \\
-c_{1,2} x_{2,0} - c_{2,2} x_{1,0} + 2\omega s_{3,2} + \dfrac{s_{1,1} s_{2,1}}{2} + \dfrac{8 c_{3,2}}{3} - \dfrac{c_{1,1} c_{2,1}}{2} = 0, \\
-s_{1,2} x_{2,0} - s_{2,2} x_{1,0} + \dfrac{8 s_{3,2}}{3} - \dfrac{c_{1,1} s_{2,1}}{2} - \dfrac{c_{2,1} s_{1,1}}{2} - 2 c_{3,2}\omega = 0, \\
\dfrac{8 x_{3,0}}{3} - x_{1,0} x_{2,0} - \dfrac{s_{1,2} s_{2,2}}{2} - \dfrac{s_{1,1} s_{2,1}}{2} - \dfrac{c_{1,2} c_{2,2}}{2} - \dfrac{c_{1,1} c_{2,1}}{2} = 0, \\
x_{3,0} + c_{3,1} + c_{3,2} - 27 = 0.
\end{cases}$$

Отметим, что для любого $h$ подобная система имеет решения

$$x_{1,0} = x_{2,0} = \pm\sqrt{b(r-1)},\ x_{3,0} = r-1,\ c_{k,i} = 0,\ s_{k,i} = 0,$$
$$\omega\ -\ \text{любое число},\ k = \overline{1,3},\ i = \overline{1,h},$$

соответствующие указанным положениям равновесия (1.3).

Таким образом, полученная нелинейная система алгебраических уравнений имеет неединственное решение. Для отыскания её приближённых решений будем использовать численный метод Ньютона, сходимость которого к нужному решению (т. е. отличному от положения равновесия) зависит от выбора начального приближения.

## 2. Символьные вычисления для получения системы алгебраических уравнений

Итак, для получения приближения к периодическому решению мы должны получить нелинейную систему относительно неизвестных коэффициентов разложения и частоты. Как показано в п. 1, даже для двух гармоник система имеет громоздкий вид. Поэтому рассмотрим алгоритм проведения символьных вычислений для её получения.

При разработке программного обеспечения был выбран математический пакет Maxima. Программа получения амплитуд и постоянных членов невязок при $h = 2$ представлена далее.

```
/* [wxMaxima batch file version 1] [ DO NOT EDIT BY HAND! ]*/
/* [wxMaxima: input   start ] */
display2d:false$
x1:x10+c1c1*cos(1*omega*t)+s1c1*sin(1*omega*t)+
c1c2*cos(2*omega*t)+s1c2*sin(2*omega*t)$
x2:x20+c2c1*cos(1*omega*t)+s2c1*sin(1*omega*t)+
c2c2*cos(2*omega*t)+s2c2*sin(2*omega*t)$
x3:x30+c3c1*cos(1*omega*t)+s3c1*sin(1*omega*t)+
c3c2*cos(2*omega*t)+s3c2*sin(2*omega*t)$
assume(omega > 0)$
delta1:trigreduce(diff(x1,t)-(10*(x2-x1)),t)$
delta2:trigreduce(diff(x2,t)-(28*x1-x2-x1*x3),t)$
delta3:trigreduce(diff(x3,t)-(x1*x2-8/3*x3),t)$
expand(diff(delta1,cos(1*omega*t)));
expand(diff(delta1,sin(1*omega*t)));
expand(diff(delta1,cos(2*omega*t)));
expand(diff(delta1,sin(2*omega*t)));
expand(integrate(delta1,t,0,2*%pi/omega)*omega/(2*%pi));
expand(diff(delta2,cos(1*omega*t)));
expand(diff(delta2,sin(1*omega*t)));
expand(diff(delta2,cos(2*omega*t)));
```

```
expand(diff(delta2,sin(2*omega*t)));
expand(integrate(delta2,t,0,2*%pi/omega)*omega/(2*%pi));
expand(diff(delta3,cos(1*omega*t)));
expand(diff(delta3,sin(1*omega*t)));
expand(diff(delta3,cos(2*omega*t)));
expand(diff(delta3,sin(2*omega*t)));
expand(integrate(delta3,t,0,2*%pi/omega)*omega/(2*%pi));
/* [wxMaxima: input   end   ] */
```

Выражение `display2d:false$` выключает многострочное рисование дробей, степеней и т. п. Знак `$` позволяет вычислить результат выражения, но не выводить на экран (вместо `;`). Функция `trigreduce(выражение,t)` свёртывает все произведения тригонометрических функций относительно переменной $t$ в комбинации сумм. Дифференцирование невязок по гармоническим функциям необходимо для получения соответствующих амплитуд. Функция `expand(выражение)` раскрывает скобки (выполняет умножение, возведение в степень, приводит подобные слагаемые).

Для нахождения постоянных членов невязок применяется интегрирование на периоде, т. е. постоянный член $k$-ой невязки равен

$$\frac{\omega \int_0^{\frac{2\pi}{\omega}} \delta_k(t)dt}{2\pi}.$$

Чтобы при символьном интегрировании пакет не задавал вопрос о знаке частоты, дается команда `assume(omega > 0)$`.

Файл с командами пакета формируется аналогично для любого количества $h$ гармоник. После выполнения данной программы пакет выведет в консоли символьные выражения для левой части системы алгебраических уравнений, которая будет решаться в нём же методом Ньютона.

Заметим, что самая затратная по времени операция здесь — символьное интегрирование. Например, для 120 гармоник время формирования системы — более 2-х суток. Здесь можно распараллелить вычислительный процесс на три компьютера, но это значительного эффекта не даст. Поэтому систему алгебраических уравнений нужно формировать сразу. Далее получим общий вид этой системы.

Отметим, что при решении системы нелинейных уравнений методом Ньютона матрица Якоби для левой части системы не обращается — в пакете Maxima используется LU-разложение для решения системы линейных уравнений на каждой итерации метода.

### 3. Общий вид системы алгебраических уравнений

Поскольку правая часть системы (0.1) содержит нелинейности в виде произведений фазовых координат, получим соотношения, выражающие коэффициенты тригонометрических полиномов, получаемых при умножении приближений $\tilde{x}_1(t)\tilde{x}_3(t)$ и $\tilde{x}_1(t)\tilde{x}_2(t)$.

Рассмотрим две функции $f(t)$ и $F(t)$, представимые рядами Фурье

$$f(t) = a_0 + \sum_{i=1}^{\infty} \left(a_i \cos(i\omega t) + b_i \sin(i\omega t)\right),$$

$$F(t) = A_0 + \sum_{i=1}^{\infty} \left(A_i \cos(i\omega t) + B_i \sin(i\omega t)\right).$$

Пусть

$$f(t)F(t) = \alpha_0 + \sum_{i=1}^{\infty} \left(\alpha_i \cos(i\omega t) + \beta_i \sin(i\omega t)\right).$$

Следуя книге [14, с. 123–125], имеем следующие соотношения:

$$\alpha_0 = a_0 A_0 + \frac{1}{2} \sum_{m=1}^{\infty} \left(a_m A_m + b_m B_m\right),$$

$$\alpha_i = a_0 A_i + \frac{1}{2} \sum_{m=1}^{\infty} \left(a_m(A_{m+i} + A_{m-i}) + b_m(B_{m+i} + B_{m-i})\right), \qquad (3.4)$$

$$\beta_i = a_0 B_i + \frac{1}{2} \sum_{m=1}^{\infty} \left(a_m(B_{m+i} - B_{m-i}) - b_m(A_{m+i} - A_{m-i})\right). \qquad (3.5)$$

Будем предполагать, что при $i > h$

$$a_i = b_i = A_i = B_i = 0.$$

Поскольку для нашей задачи мы ищем приближение до $h$-ой гармоники включительно, занулим все амплитуды в произведении при $i > h$, т. е.

$$\alpha_i = \beta_i = 0.$$

Таким образом, мы перейдём от произведения рядов к произведению тригонометрических полиномов.

Также в соотношениях (3.4) и (3.5) мы будем предполагать [14, с. 124], что

$$A_{m-i} = A_{i-m}, \ B_{m-i} = -B_{i-m}, \ B_0 = 0.$$

Тогда получим:

$$\alpha_0 = a_0 A_0 + \frac{1}{2} \sum_{m=1}^{h} \left(a_m A_m + b_m B_m\right),$$

$$\alpha_i = a_0 A_i + \frac{1}{2}\sum_{m=1}^{\infty} a_m A_{m+i} + \frac{1}{2}\sum_{m=1}^{\infty} a_m A_{m-i} + \frac{1}{2}\sum_{m=1}^{\infty} b_m B_{m+i} + \frac{1}{2}\sum_{m=1}^{\infty} b_m B_{m-i} =$$

$$= a_0 A_i + \frac{1}{2}\sum_{m=1}^{h-i} a_m A_{m+i} + \frac{1}{2} a_i A_0 + \frac{1}{2}\sum_{m=1}^{i-1} a_m A_{i-m} + \frac{1}{2}\sum_{m=i+1}^{h} a_m A_{m-i} +$$

$$+ \frac{1}{2}\sum_{m=1}^{h} b_m B_{m+i} + \frac{1}{2} b_i B_0 - \frac{1}{2}\sum_{m=1}^{i-1} b_m B_{i-m} + \frac{1}{2}\sum_{m=i+1}^{h} b_m B_{m-i} =$$

$$= a_0 A_i + a_i A_0 + \frac{1}{2}\sum_{m=1}^{h-i}(a_m A_{m+i} + b_m B_{m+i}) + \frac{1}{2}\sum_{m=1}^{i-1}(a_m A_{i-m} - b_m B_{i-m}) +$$

$$+ \frac{1}{2}\sum_{m=i+1}^{h}(a_m A_{m-i} + b_m B_{m-i}),$$

$$\beta_i = a_0 B_i + \frac{1}{2}\sum_{m=1}^{\infty} a_m B_{m+i} - \frac{1}{2}\sum_{m=1}^{\infty} a_m B_{m-i} - \frac{1}{2}\sum_{m=1}^{\infty} b_m A_{m+i} + \frac{1}{2}\sum_{m=1}^{\infty} b_m A_{m-i} =$$

$$= a_0 B_i + \frac{1}{2}\sum_{m=1}^{h-i} a_m B_{m+i} + \frac{1}{2}\sum_{m=1}^{i-1} a_m B_{i-m} - \frac{1}{2}\sum_{m=i+1}^{h} a_m B_{m-i} -$$

$$- \frac{1}{2}\sum_{m=1}^{h-i} b_m A_{m+i} + b_i A_0 + \frac{1}{2}\sum_{m=1}^{i-1} b_m A_{i-m} + \frac{1}{2}\sum_{m=i+1}^{h} b_m A_{m-i} =$$

$$= a_0 B_i + b_i A_0 + \frac{1}{2}\sum_{m=1}^{h-i}(a_m B_{m+i} - b_m A_{m+i}) + \frac{1}{2}\sum_{m=1}^{i-1}(a_m B_{i-m} + b_m A_{i-m}) +$$

$$+ \frac{1}{2}\sum_{m=i+1}^{h}(-a_m B_{m-i} + b_m A_{m-i}).$$

Применяя полученные формулы для вычисления произведений тригонометрических полиномов к невязкам, мы можем записать уравнения для $i$-ых гармоник ($i = \overline{1,h}$ — номер гармоники, $k = \overline{1,3}$ — номер невязки):

$k = 1$:

$$i\omega s_{1,i} - 10 c_{2,i} + 10 c_{1,i} = 0,$$
$$-i\omega c_{1,i} - 10 s_{2,i} + 10 s_{1,i} = 0,$$

уравнение, соответствующее постоянному члену для первой невязки, —

$$x_{1,0} - x_{2,0} = 0,$$

$k=2$:

$$i\omega s_{2,i} - 28c_{1,i} + c_{2,i} + x_{1,0}c_{3,i} + c_{1,i}x_{3,0} + \frac{1}{2}\sum_{m=1}^{h-i}(c_{1,m}c_{3,m+i} + s_{1,m}s_{3,m+i}) +$$
$$+ \frac{1}{2}\sum_{m=1}^{i-1}(c_{1,m}c_{3,i-m} - s_{1,m}s_{3,i-m}) +$$
$$+ \frac{1}{2}\sum_{m=i+1}^{h}(c_{1,m}c_{3,m-i} + s_{1,m}s_{3,m-i}) = 0,$$

$$-i\omega c_{2,i} - 28s_{1,i} + s_{2,i} + x_{1,0}s_{3,i} + s_{1,i}x_{3,0} + \frac{1}{2}\sum_{m=1}^{h-i}(c_{1,m}s_{3,m+i} - s_{1,m}c_{3,m+i}) +$$
$$+ \frac{1}{2}\sum_{m=1}^{i-1}(c_{1,m}s_{3,i-m} + s_{1,m}c_{3,i-m}) +$$
$$+ \frac{1}{2}\sum_{m=i+1}^{h}(-c_{1,m}s_{3,m-i} + s_{1,m}c_{3,m-i}) = 0,$$

уравнение, соответствующее постоянному члену для второй невязки, —

$$-28x_{1,0} + x_{2,0} + x_{1,0}x_{3,0} + \frac{1}{2}\sum_{m=1}^{h}(c_{1,m}c_{3,m} + s_{1,m}s_{3,m}) = 0,$$

$k=3$:

$$i\omega s_{3,i} - x_{1,0}c_{2,i} - c_{1,i}x_{2,0} - \frac{1}{2}\sum_{m=1}^{h-i}(c_{1,m}c_{2,m+i} + s_{1,m}s_{2,m+i}) -$$
$$- \frac{1}{2}\sum_{m=1}^{i-1}(c_{1,m}c_{2,i-m} - s_{1,m}s_{2,i-m}) -$$
$$- \frac{1}{2}\sum_{m=i+1}^{h}(c_{1,m}c_{2,m-i} + s_{1,m}s_{2,m-i}) + \frac{8}{3}c_{3,i} = 0,$$

$$-i\omega c_{3,i} - x_{1,0}s_{2,i} - s_{1,i}x_{2,0} - \frac{1}{2}\sum_{m=1}^{h-i}(c_{1,m}s_{2,m+i} - s_{1,m}c_{2,m+i}) -$$
$$- \frac{1}{2}\sum_{m=1}^{i-1}(c_{1,m}s_{2,i-m} + s_{1,m}c_{2,i-m}) -$$
$$- \frac{1}{2}\sum_{m=i+1}^{h}(-c_{1,m}s_{2,m-i} + s_{1,m}c_{2,m-i}) + \frac{8}{3}s_{3,i} = 0,$$

уравнение, соответствующее постоянному члену для третьей невязки, —

$$-x_{1,0}x_{2,0} - \frac{1}{2}\sum_{m=1}^{h}(c_{1,m}c_{2,m} + s_{1,m}s_{2,m}) + \frac{8}{3}x_{3,0} = 0,$$

дополнительное уравнение системы —

$$x_{3,0} + \sum_{i=1}^{h} c_{3,i} - 27 = 0.$$

## 4. Результаты вычислительного эксперимента

В результате многочисленных вычислительных экспериментов было подобрано начальное приближение для циклической частоты, постоянных членов и амплитуд при $h = h_1 = 5$:

$$\omega = 4, \ x_{1,0} = x_{2,0} = x_{3,0} = 0, \ c_{1,i} = -1, \ i = \overline{1,5},$$
$$s_{1,j} = 0, \ j = 1,3,4,5, \ s_{1,2} = 1.$$

Данный результат замечателен тем, что метод Ньютона сходится к решению, отличному от положений равновесия. Поэтому для улучшения точности приближенного периодического решения мы рассматриваем систему алгебраических уравнений для значения $h$, равного некоторому $h_2 > h_1$. Полученное численное решение системы при $h = h_1$ берётся как начальное приближение для амплитуд с индексами $i \leq h_1$ у системы с $h = h_2$, а значения начального приближения для амплитуд с индексами $i > h_1$ полагаются равными нулю.

В таблицах 1–3 приведён результат решения системы при $h = 35$, точность метода Ньютона – $10^{-8}$. Значение периода получается равным $T = 1.558652210$, начальное условие для полученного приближённого периодического решения –

$$\tilde{x}_1(0) = -2.147367631, \ \tilde{x}_2(0) = 2.078048211, \ \tilde{x}_3(0) = 27. \tag{4.6}$$

Начальные значения (4.6) были проверены на периоде в компьютерной программе, реализующей численное интегрирование системы (0.1) модифицированным методом степенных рядов [9] с точностью оценки общего члена ряда $10^{-25}$, 100 бит под мантиссу вещественного числа и машинным эпсилон $1.57772 \cdot 10^{-30}$. При таких параметрах метода приближённые значения фазовых координат, полученные с помощью численного интегрирования, были также проверены тем же численным методом, но в обратном времени. Значения в обратном времени совпадают с (4.6) до 9-го знака включительно после точки. Результирующие же значения $x_1(T)$, $x_2(T)$ и $x_3(T)$ совпадают с (4.6) до 8-го знака включительно.

Таблица 1

Амплитуды гармоник для $\tilde{x}_1(t)$, $x_{1,0} = 0$

| $i$ | $c_{1,i}$ | $s_{1,i}$ |
|---|---|---|
| 1 | $-5.780478259196228$ | $8.56017654325353$ |
| 2 | $0$ | $0$ |
| 3 | $3.160762628380509$ | $2.239212141102876$ |
| 4 | $0$ | $0$ |
| 5 | $0.6958870387616096$ | $-0.7979388979225431$ |
| 6 | $0$ | $0$ |
| 7 | $-0.1891992374027477$ | $-0.1864921358925765$ |
| 8 | $0$ | $0$ |
| 9 | $-0.04770429623010056$ | $0.04554044367245914$ |
| 10 | $0$ | $0$ |
| 11 | $0.01112322884679491$ | $0.01209138588669679$ |
| 12 | $0$ | $0$ |
| 13 | $0.003061207095371694$ | $-0.002735092350544739$ |
| 14 | $0$ | $0$ |
| 15 | $-6.744578887916229 \cdot 10^{-4}$ | $-7.748319471034087 \cdot 10^{-4}$ |
| 16 | $0$ | $0$ |
| 17 | $-1.960718247379475 \cdot 10^{-4}$ | $1.665584161919807 \cdot 10^{-4}$ |
| 18 | $0$ | $0$ |
| 19 | $4.116738805347028 \cdot 10^{-5}$ | $4.960493476144467 \cdot 10^{-5}$ |
| 20 | $0$ | $0$ |
| 21 | $1.254757391175977 \cdot 10^{-5}$ | $-1.018054283421179 \cdot 10^{-5}$ |
| 22 | $0$ | $0$ |
| 23 | $-2.518375902000733 \cdot 10^{-6}$ | $-3.173486439630506 \cdot 10^{-6}$ |
| 24 | $0$ | $0$ |
| 25 | $-8.025338211960923 \cdot 10^{-7}$ | $6.230623750431923 \cdot 10^{-7}$ |
| 26 | $0$ | $0$ |
| 27 | $1.541534734542893 \cdot 10^{-7}$ | $2.0292802821633 \cdot 10^{-7}$ |
| 28 | $0$ | $0$ |
| 29 | $5.130649139299358 \cdot 10^{-8}$ | $-3.813725452268523 \cdot 10^{-8}$ |
| 30 | $0$ | $0$ |
| 31 | $-9.43393531993558 \cdot 10^{-9}$ | $-1.297038481588497 \cdot 10^{-8}$ |
| 32 | $0$ | $0$ |
| 33 | $-3.278552746800046 \cdot 10^{-9}$ | $2.333260259021725 \cdot 10^{-9}$ |
| 34 | $0$ | $0$ |
| 35 | $5.76957885768651 \cdot 10^{-10}$ | $8.28626640138045 \cdot 10^{-10}$ |

Таблица 2

Амплитуды гармоник для $\tilde{x}_2(t)$, $x_{2,0} = 0$

| $i$ | $c_{2,i}$ | $s_{2,i}$ |
|---|---|---|
| 1 | $-2.32972926505593$ | $10.89038310357172$ |
| 2 | 0 | 0 |
| 3 | $5.86875317198698$ | $-1.5832552129833$ |
| 4 | 0 | 0 |
| 5 | $-0.9124249133801483$ | $-2.200556873678218$ |
| 6 | 0 | 0 |
| 7 | $-0.7154457265566421$ | $0.3473932955614448$ |
| 8 | 0 | 0 |
| 9 | $0.1175186702136983$ | $0.2186139734768588$ |
| 10 | 0 | 0 |
| 11 | $0.06473984670858603$ | $-0.03723215039412078$ |
| 12 | 0 | 0 |
| 13 | $-0.01127208646321726$ | $-0.01877739524860192$ |
| 14 | 0 | 0 |
| 15 | $-0.005359671824365359$ | $0.003303445299126894$ |
| 16 | 0 | 0 |
| 17 | $9.453499475830811 \cdot 10^{-4}$ | $0.001510235036151227$ |
| 18 | 0 | 0 |
| 19 | $4.211022386354685 \cdot 10^{-4}$ | $-2.657049331814368 \cdot 10^{-4}$ |
| 20 | 0 | 0 |
| 21 | $-7.363528144366622 \cdot 10^{-5}$ | $-1.164013765469982 \cdot 10^{-4}$ |
| 22 | 0 | 0 |
| 23 | $-3.19419300699788 \cdot 10^{-5}$ | $2.017609175377016 \cdot 10^{-5}$ |
| 24 | 0 | 0 |
| 25 | $5.47663534401654 \cdot 10^{-6}$ | $8.710929378319451 \cdot 10^{-6}$ |
| 26 | 0 | 0 |
| 27 | $2.362852034076972 \cdot 10^{-6}$ | $-1.474901091428546 \cdot 10^{-6}$ |
| 28 | 0 | 0 |
| 29 | $-3.94532524722541 \cdot 10^{-7}$ | $-6.379296603810031 \cdot 10^{-7}$ |
| 30 | 0 | 0 |
| 31 | $-1.715198229248314 \cdot 10^{-7}$ | $1.049218598356554 \cdot 10^{-7}$ |
| 32 | 0 | 0 |
| 33 | $2.776045093375681 \cdot 10^{-8}$ | $4.59473450493284 \cdot 10^{-8}$ |
| 34 | 0 | 0 |
| 35 | $1.22681173575872 \cdot 10^{-8}$ | $-7.31171826830086 \cdot 10^{-9}$ |

Таблица 3

Амплитуды гармоник для $\tilde{x}_3(t)$, $x_{3,0} = 23.04210397942006$

| $i$ | $c_{3,i}$ | $s_{3,i}$ |
|---|---|---|
| 1 | 0 | 0 |
| 2 | 7.568410271550653 | $-9.50386584559212$ |
| 3 | 0 | 0 |
| 4 | $-3.555327211552558$ | $-1.844710563805469$ |
| 5 | 0 | 0 |
| 6 | $-0.4741220131932616$ | 1.279043179069961 |
| 7 | 0 | 0 |
| 8 | 0.4227292179138024 | 0.1274574086305204 |
| 9 | 0 | 0 |
| 10 | 0.03498415351761577 | $-0.1315337800809524$ |
| 11 | 0 | 0 |
| 12 | $-0.03934013541135439$ | $-0.009645786231708874$ |
| 13 | 0 | 0 |
| 14 | $-0.002660522588813564$ | 0.01145537653603837 |
| 15 | 0 | 0 |
| 16 | 0.003271688724557337 | $7.33752523103949 \cdot 10^{-4}$ |
| 17 | 0 | 0 |
| 18 | $2.024982256871223 \cdot 10^{-4}$ | $-9.206266886554897 \cdot 10^{-4}$ |
| 19 | 0 | 0 |
| 20 | $-2.560063570343799 \cdot 10^{-4}$ | $-5.58964460662525 \cdot 10^{-5}$ |
| 21 | 0 | 0 |
| 22 | $-1.542436654918173 \cdot 10^{-5}$ | $7.050327849098175 \cdot 10^{-5}$ |
| 23 | 0 | 0 |
| 24 | $1.926014222030195 \cdot 10^{-5}$ | $4.25261452471065 \cdot 10^{-6}$ |
| 25 | 0 | 0 |
| 26 | $1.170939944189529 \cdot 10^{-6}$ | $-5.225643926851625 \cdot 10^{-6}$ |
| 27 | 0 | 0 |
| 28 | $-1.409525591131397 \cdot 10^{-6}$ | $-3.21879984959824 \cdot 10^{-7}$ |
| 29 | 0 | 0 |
| 30 | $-8.83134288999026 \cdot 10^{-8}$ | $3.782652721710986 \cdot 10^{-7}$ |
| 31 | 0 | 0 |
| 32 | $1.010610960272394 \cdot 10^{-7}$ | $2.418021923473667 \cdot 10^{-8}$ |
| 33 | 0 | 0 |
| 34 | $6.606163280924149 \cdot 10^{-9}$ | $-2.689431432873997 \cdot 10^{-8}$ |
| 35 | 0 | 0 |

Цикл, соответствующий (4.6), показан на рис. 1. Отметим, что найденный цикл совпадает с первым циклом Вишваната в [7], все знаки после точки для величины $T$ также совпадают с данными из работы [7].

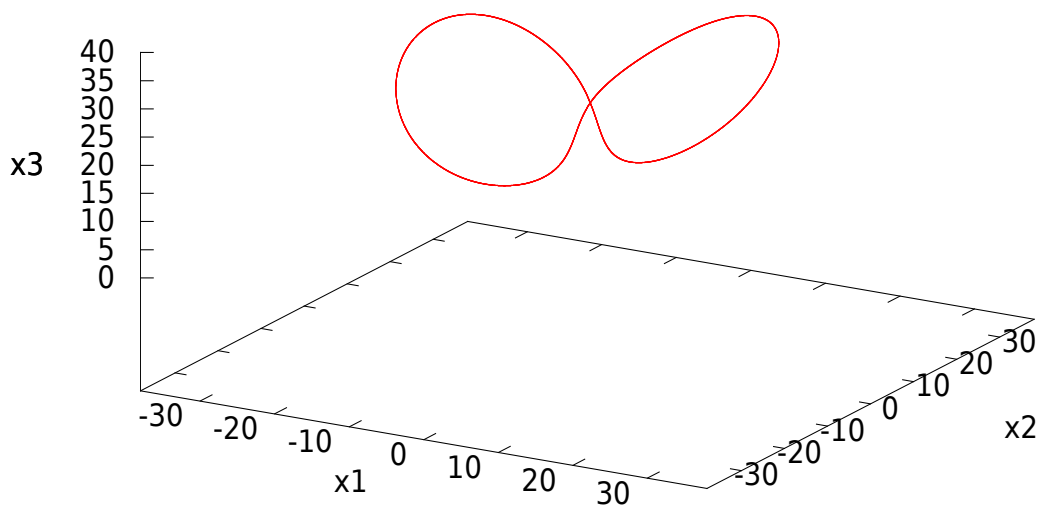

**Рис. 1.** Цикл, полученный методом гармонического баланса

## Список литературы

**Информация об авторах**

**Пчелинцев Александр Николаевич**, кандидат физико-математических наук, заведующий кафедрой «Высшая математика». Тамбовский государственный технический университет, г. Тамбов, Российская Федерация. E-mail: pchelintsev.an@yandex.ru
**ORCID:** https://orcid.org/0000-0003-4136-1227

**Полуновский Андрей Андреевич**, студент. Московский государственный технический университет им. Н. Э. Баумана (национальный исследовательский университет), г. Москва, Российская Федерация. E-mail: apap2009@yandex.ru
**ORCID:** https://orcid.org/0000-0002-6557-3649

**Юханова Ирина Юрьевна**, магистр. Тамбовский государственный технический университет, г. Тамбов, Российская Федерация. E-mail: irina_yu_10@mail.ru
**ORCID:** https://orcid.org/0000-0002-8339-0459

**Information about the authors**

**Alexander N. Pchelintsev**, Candidate of Physics and Mathematics, Head of the Higher Mathematics Department. Tambov State Technical University, Tambov, the Russian Federation. E-mail: pchelintsev.an@yandex.ru
**ORCID:** https://orcid.org/0000-0003-4136-1227

**Andrey A. Polunovskiy**, Student. Bauman Moscow State Technical University, Moscow, the Russian Federation. E-mail: apap2009@yandex.ru
**ORCID:** https://orcid.org/0000-0002-6557-3649

**Irina Yu. Yukhanova**, Master. Tambov State Technical University, Tambov, the Russian Federation. E-mail: irina_yu_10@mail.ru
**ORCID:** https://orcid.org/0000-0002-8339-0459